\documentclass[12pt,a4paper]{article}
\usepackage[top=20mm, left=25mm, right=25mm, bottom=20mm]{geometry}
\usepackage{graphicx} 
\usepackage[colorlinks = true]{hyperref}
\usepackage[utf8]{inputenc}
\usepackage{xcolor}
\usepackage{bm}
\usepackage{amsmath}
\usepackage{amsfonts}
\usepackage{amssymb}
\usepackage{amsthm}
\usepackage{tikz}
\usetikzlibrary{shapes,arrows,positioning,calc}

\usepackage{svg} 

\newtheorem{theorem}{Theorem}[section]

\newtheorem{corollary}[theorem]{Corollary}

\theoremstyle{definition}
\newtheorem{definition}[theorem]{Definition}
\newtheorem{example}[theorem]{Example}
\newtheorem{remark}[theorem]{Remark}

\usepackage{enumitem}
\usepackage{mathtools}

\newtheorem{observation}[theorem]{Observation}

\title{Stronger Spectral Bound for Clique Number}

\author{
  Hareshkumar Jadav\textsuperscript{1}\thanks{Corresponding author. Email: \texttt{harershjadavcse@gmail.com}}
  \and
  Sudipta Mallik\textsuperscript{2}\thanks{Email: \texttt{mallik@marshall.edu}}
  \and
  Akhilesh Panchangam\textsuperscript{1}
  \and
  Ranveer Singh\textsuperscript{1}\thanks{Email: \texttt{ranveer@iiti.ac.in}}
  \\
  \textsuperscript{1}Department of Computer Science and Engineering, IIT Indore, India \\
  \textsuperscript{2}Department of Mathematics and Physics, Marshall University, USA
}

\date{}

\begin{document}

\maketitle
\begin{abstract}
    Given an integer $k$, deciding whether a graph has a clique of size $k$ is an NP-complete problem. Wilf's inequality provides a spectral lower bound for the clique number of a graph. Strengthening this bound, Elphick and Wocjan proposed a conjecture in 2018 using positive adjacency eigenvalues. In this paper, we present an infinite family of graphs  for which the Elphick and Wocjan conjecture is true. We find a stronger spectral lower bound for the clique number for this family of graphs. 
\end{abstract}

2020 Mathematics Subject Classification:  05C50 

Keywords: Clique number, Eigenvalue, Graph

\section{Introduction}
Let $G=(V,E)$ be a simple graph with $n=|V|$ vertices and $m=|E|$ edges. Let $A(G)$ denote the adjacency matrix of $G$ and write its eigenvalues as $\lambda_1\ge \lambda_2\ge \cdots \ge \lambda_n$. A \emph{clique} is a set of pairwise adjacent vertices, i.e., a complete graph. The size of a maximum clique in $G$ is the \emph{clique number} of $G$, denoted by $\omega(G)$. The decision problem “does $G$ contain a clique of size at least $k$?” is NP-complete. Spectral bounds provide efficiently computable lower bounds for $\omega$. 
In particular, Wilf~\cite{wilf} proved that 
\begin{equation}\label{eq:wilf-lambda}
\lambda_{1} \;\le\; n\!\left(1-\frac{1}{\omega}\right).
\end{equation}
Equivalently,
\begin{equation}\label{eq:wilf-omega}
\omega \;\ge\; \frac{n}{\,n-\lambda_{1}\,}.
\end{equation}
Elphick and coauthors observed that many classical eigenvalue bounds for the chromatic and clique numbers can be sharpened by replacing the term $\lambda_1$ 
with 
\[
\sqrt{s^{+}} \;=\; \sqrt{\sum_{\lambda_i>0}\lambda_i^{2}}.
\]

where the sum runs over all positive eigenvalues (counted with multiplicity).
For example, Edwards and Elphick~\cite{edwards1983lower} gave a lower bound on $\chi(G)$ in terms of $\lambda_1$ and Ando and Lin~\cite{ando2015proof} later proved a conjecture of Wocjan and Elphick~\cite{wocjan2012new} showing that the same type of bound follow with $\sqrt{s^{+}}$ in place of $\lambda_1$. 
Further results of Elphick et al.\ and Wu and Elphick show that several other inequalities involving $\chi(G)$, $\omega(G)$ and $\lambda_1$ admit stronger versions when $\lambda_1$ is systematically replaced by $\sqrt{s^{+}}$; see, e.g.,~\cite{favaron1993some,stanley1987bound,wu2017upper}.

Motivated by this “$\lambda_1 \to \sqrt{s^{+}}$’’ philosophy, Elphick and Wocjan \cite{elphick2024two} proposed the following strengthening of Wilf’s bound, which we refer to as the \emph{Elphick–Wocjan conjecture}:
\begin{equation}\label{conj2}
\sqrt{s^{+}} \;\le\; n\!\left(1-\frac{1}{\omega}\right)
\end{equation}
Equivalently,
\begin{equation}\label{conj2-omega}
\omega \;\ge\; \frac{n}{\,n-\sqrt{s^{+}}\,}.
\end{equation}
They provided extensive computational evidence on thousands of named graphs up to 40 vertices and proved the conjecture for several graph classes; see \cite{elphick2024two}. Independently, \cite{jadav2025strengthening} verifies the Elphick--Wocjan inequality for additional families, including conference graphs, some strongly regular graphs, line graphs of $K_n$ and certain Cartesian squares.


In this work, our main result is a strengthened spectral lower bound on the clique number, in the spirit of Wilf’s inequality. Specifically, for two infinite families of graphs, denoted by $B_3(n_1, n_2, n_3)$ and $B_4(n_1, n_2, n_3, n_4)$ (formally defined in Section~\ref{section2}), we show that the term $\lambda_1$ in Wilf’s bound can be
replaced by
\[
  \sqrt{\sum_{\lambda_i \neq -1} \lambda_i^{2}},
\]
the sum of squares of all adjacency eigenvalues except a single negative eigenvalue $-1$. This replacement yields a sharper spectral lower bound on $\omega(G)$ for every graph in these two families. Both families are defined for arbitrary positive integers $n_i$ and consist of graphs with exactly two positive eigenvalues.
As a corollary, we obtain that the Elphick--Wocjan conjecture holds for all graphs in these families. This paper is organized as follows. 
In Section~\ref{section2}, we introduce the graph families and establish several supporting results needed for our main theorem. 
In Section~\ref{sec:main_result}, we prove the main result and show that, as a consequence, Elphick–Wocjan conjecture also holds for this class of graphs.



\section{Graphs Structure and Notations}\label{section2}
The infinite graph class we consider is closely related to graphs with exactly two positive eigenvalues and zero nullity which are studied in \cite{Oboudi}. We first introduce some notations and see the structure of graphs in this graph class. 

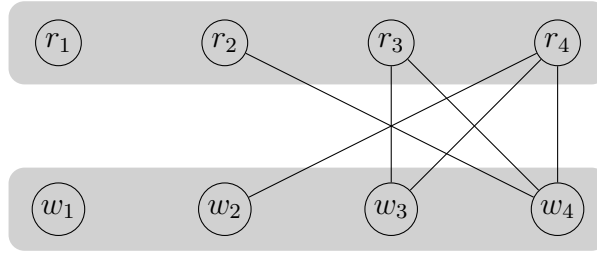
\begin{figure}
    \centering

\begin{tikzpicture}[scale=1.1,
  every node/.style={circle, draw, inner sep=1.8pt}
]

\fill[black!30, opacity=0.6, rounded corners=6pt]
      (-0.6,1.5) rectangle (6.6,2.5);

\node (v1) at (0,2) {$r_1$};
\node (v2) at (2,2) {$r_2$};
\node (v3) at (4,2) {$r_3$};
\node (v4) at (6,2) {$r_4$};

\fill[black!30, opacity=0.6, rounded corners=6pt]
      (-0.6,-0.5) rectangle (6.6,0.5);

\node (w1) at (0,0) {$w_1$};
\node (w2) at (2,0) {$w_2$};
\node (w3) at (4,0) {$w_3$};
\node (w4) at (6,0) {$w_4$};

\draw (v2) -- (w4);

\draw (v3) -- (w4);
\draw (v3) -- (w3);

\draw (v4) -- (w4);
\draw (v4) -- (w3);
\draw (v4) -- (w2);

\end{tikzpicture}

    \caption{The graph $G_8$. 
    The gray shading indicates that the vertices inside each block induce a complete graph.}
    \label{fig:G8}
\end{figure}

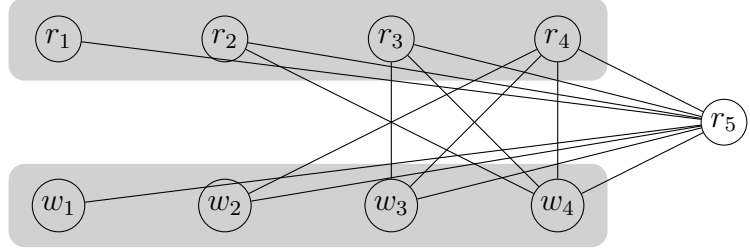
\begin{figure}
    \centering

\begin{tikzpicture}[scale=1.1,
  every node/.style={circle, draw, inner sep=1.8pt}
]

\fill[black!30, opacity=0.6, rounded corners=6pt]
      (-0.6,1.5) rectangle (6.6,2.5);

\node (v1) at (0,2) {$r_1$};
\node (v2) at (2,2) {$r_2$};
\node (v3) at (4,2) {$r_3$};
\node (v4) at (6,2) {$r_4$};
\node (v5) at (8,1) {$r_5$};

\fill[black!30, opacity=0.6, rounded corners=6pt]
      (-0.6,-0.5) rectangle (6.6,0.5);

\node (w1) at (0,0) {$w_1$};
\node (w2) at (2,0) {$w_2$};
\node (w3) at (4,0) {$w_3$};
\node (w4) at (6,0) {$w_4$};

\draw (v2) -- (w4);

\draw (v3) -- (w4);
\draw (v3) -- (w3);

\draw (v4) -- (w4);
\draw (v4) -- (w3);
\draw (v4) -- (w2);

\draw (v5) -- (w4);
\draw (v5) -- (w3);
\draw (v5) -- (w2);
\draw (v5) -- (w1);
\draw (v5) -- (v4);
\draw (v5) -- (v3);
\draw (v5) -- (v2);
\draw (v5) -- (v1);

\end{tikzpicture}

    \caption{The graph $G_9$. 
    The gray shading indicates that the vertices inside each block induce a complete graph.}
    \label{fig:G9}
\end{figure}

\begin{definition}\cite{Duan23}\label{def_Gn}
For an integer $t \geq 2$, let $K_{\lceil \frac{t}{2} \rceil}$ and $K_{\lfloor \frac{t}{2} \rfloor}$ be vertex-disjoint complete graphs with vertex sets
\[
V = \{r_1, r_2, \ldots, r_{\lceil \frac{t}{2} \rceil}\}
\quad\text{and}\quad
W = \{w_1, w_2, \ldots, w_{\lfloor \frac{t}{2} \rfloor}\},
\]
respectively. Let $N_W(x)$ denote the neighborhood of any vertex $x$ in $W$. Let  $G_t$ (see Figure \ref{fig:G8} for $t=8$ and Figure \ref{fig:G9} for $t=9$) be the graph obtained by applying the following operations:
\begin{itemize}
\item If $t$ is even, we insert new edges such that
\begin{align*}
&N_W(r_1) = \varnothing \subset N_W(r_2) = \{w_{\frac{t}{2}}\} \subset N_W(r_3) = \{w_{\frac{t}{2}}, w_{\frac{t}{2}-1}\} \subset \cdots \\
&\subset N_W(r_{\frac{t}{2}-1}) = \{w_{\frac{t}{2}}, w_{\frac{t}{2}-1}, \ldots, w_3\} 
\subset N_W(r_{\frac{t}{2}}) = \{w_{\frac{t}{2}}, w_{\frac{t}{2}-1}, \ldots, w_2\}.
\end{align*}
\item For $t$ odd, we insert new edges such that
\begin{align*}
&N_W(r_1) = \varnothing \subset N_W(r_2) = \{w_{\frac{t-1}{2}}\} \subset N_W(r_3) = \{w_{\frac{t-1}{2}}, w_{\frac{t-1}{2}-1}\} \subset \cdots \\
&\subset N_W(r_{\frac{t+1}{2}-1}) = \{w_{\frac{t-1}{2}}, w_{\frac{t-1}{2}-1}, \ldots, w_2\} 
\subset N_W(r_{\frac{t+1}{2}}) = \{w_{\frac{t-1}{2}}, w_{\frac{t-1}{2}-1}, \ldots, w_1\}.
\end{align*}
\end{itemize}
Note that when $t$ is odd, the above construction of $G_t$ can equivalently be obtained using $G_{t-1}$ by adding the vertex $r_{\frac{t+1}{2}}$ and joining it to all the vertices of $G_{t-1}$. 
\end{definition}


Throughout, we fix the vertex ordering:
$V(G_{2s})=(r_1,\dots,r_s,w_1,\dots,w_s)$ for $t=2s$, and
$V(G_{2s+1})=(r_1,\dots,r_s,w_1,\dots,w_s,r_{s+1})$ for $t=2s+1$. This convention ensures that the matrix $M$ in Theorem~\ref{thm:char_polynomial} admits a nice block structure.

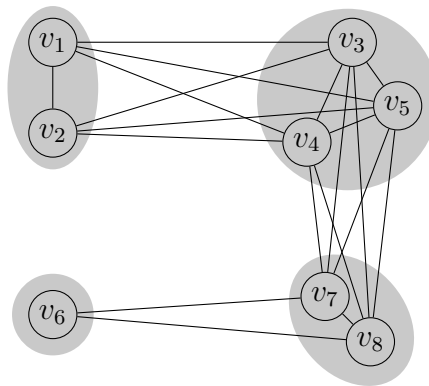
\begin{figure}[ht]
    \centering
    \begin{tikzpicture}[scale=1.2]

\tikzset{
  vertex/.style={circle, draw, inner sep=2pt}
}
\fill[black!30, opacity=0.7] (-3,4.5) ellipse (0.5cm and 0.9cm);

\fill[black!30, opacity=0.7] (0.25,4.37) ellipse (1cm and 1cm);

\fill[black!30, opacity=0.7] (-3,2) ellipse (0.45cm and 0.45cm);

\fill[
  black!30,
  opacity=0.7,
  rotate around={40:(0.15,1.95)} 
] (0.25,1.85) ellipse (0.6cm and 0.8cm);

\node[vertex] (v1) at (-3,5) {$v_1$};
\node[vertex] (v2) at (-3,4) {$v_2$};
\node[vertex] (v3) at (0.3,5) {$v_3$};
\node[vertex] (v4) at (-0.2,3.9) {$v_4$};
\node[vertex] (v5) at (0.8,4.3) {$v_5$};
\node[vertex] (v6) at (-3,2) {$v_6$};
\node[vertex] (v7) at (0,2.2) {$v_7$};
\node[vertex] (v8) at (0.5,1.7) {$v_8$};

\draw (v1) -- (v3);
\draw (v1) -- (v4);
\draw (v1) -- (v5);

\draw (v3) -- (v5);
\draw (v4) -- (v3);
\draw (v4) -- (v5);
\draw (v3) -- (v2);
\draw (v5) -- (v2);
\draw (v2) -- (v1);

\draw (v6) -- (v7);
\draw (v6) -- (v8);

\draw (v3) -- (v7);
\draw (v3) -- (v8);

\draw (v4) -- (v7);
\draw (v4) -- (v8);

\draw (v5) -- (v7);
\draw (v5) -- (v8);

\draw (v4) -- (v2);
\draw (v7) -- (v8);

\end{tikzpicture}
    \caption{The graph $B_4(2, 3, 1, 2)$. }
    \label{fig:graph-clouds}
\end{figure}

Let $G_t$ be a graph with vertex set $V(G_t)$. For positive integers
$n_1, n_2, \ldots, n_t$, we denote by
\[
    B_t({n_1}, {n_2}, \dots, {n_t})
\]
the generalized lexicographic product of $G_t$ with $K_{n_1}, K_{n_2}, \ldots, K_{n_t}$, defined as follows. For each $v_i \in V(G_t)$, we replace vertex $v_i$ by a clique $K_{n_i}$
and we make every vertex of $K_{n_i}$ adjacent to every vertex of $K_{n_j}$
whenever $v_i$ is adjacent to $v_j$ in $G_t$ (see Figure~\ref{fig:graph-clouds} for an example).

The following theorem shows that any connected graph with exactly two positive eigenvalues with zero nullity is of the form $B_t({n_1}, {n_2}, \ldots, {n_t})$.


\begin{theorem}\cite[Theorem 2.18]{Oboudi}\label{thm:Oboudi}
Let $G$ be a graph of order $n \geq 3$ with two positive eigenvalues and zero nullity. 
\begin{itemize}
    \item[(i)] If $G$ is disconnected, then $G$ is isomorphic to the disjoint union of $K_{n_1}$ and $K_{n_2}$ for some integers $n_1, n_2 \geq 2$;
    \item[(ii)] If $G$ is connected, then $G$ is isomorphic to $B_t({n_1}, {n_2}, \ldots, {n_t})$, where  $3  \leq t \leq 12$ and $n_1 + n_2 + \cdots + n_t = n$,  for some positive integers $n_1, n_2, \ldots, n_t$.
\end{itemize}
\end{theorem}

The following theorem, shows that $-1$ is always an eigenvalue of $B_t({n_1}, \ldots, {n_t})$ with multiplicity $n - t$ and expresses the remaining eigenvalues in terms of a suitable $t \times t$ matrix $M$. This will be our starting point for deriving a closed form of $\sum_{\lambda_i\ne -1}\lambda_i^2$.


\begin{theorem}\cite[Theorem 2.21]{Oboudi}\label{thm:char_polynomial} 
Let $t \ge 2$ be an integer and let $G_t$ be a graph on $t$ vertices with adjacency matrix $A = [a_{ij}]$.  Let $n_1, \ldots, n_t$ be some positive integers such that $n_1 + \cdots + n_t=n$.
Define the $t \times t$ matrix $M = [m_{ij}]$ by
\begin{equation*}
m_{ij} := 
\begin{cases}
n_i - 1, & \text{if } i = j,\\[2pt]
a_{ij}n_j, & \text{if } i \neq j.
\end{cases}
\end{equation*}
Then the eigenvalues of $B_t({n_1}, \ldots, {n_t})$  are $-1$ with multiplicity  $n - t$ and the rest of the eigenvalues are the eigenvalue of $M$. Moreover, for a fixed vertex ordering of $G_t$, the matrix $M$ admits a nice block structure.
\end{theorem}

\begin{example}\label{example1}
Consider the graph $B_4(2,3,1,2)$  with vertex set $\{v_1,\ldots,v_8\}$. We order the vertices as
\[
v_1,v_2 \in K_2,\quad
v_3,v_4,v_5 \in K_3,\quad
v_6 \in K_1,\quad
v_7,v_8 \in K_2,
\]
as illustrated in Figure~\ref{fig:graph-clouds}. In this case, $t=4$ and $(n_1,n_2,n_3,n_4) = (2,3,1,2)$.
Therefore, the matrix $M$ by  Theorem~\ref{thm:char_polynomial} is
\[
M =
\begin{bmatrix}
n_1-1 & n_2    & 0      & 0\\
n_1   & n_2-1  & 0      & n_4\\
0     & 0      & n_3-1  & n_4\\
0     & n_2    & n_3    & n_4-1
\end{bmatrix}
=
\begin{bmatrix}
1 & 3 & 0 & 0\\
2 & 2 & 0 & 2\\
0 & 0 & 0 & 2\\
0 & 3 & 1 & 1
\end{bmatrix}.
\]

The eigenvalues of $M$ are
\[
\lambda_1 \approx 5.0899,\quad
\lambda_2 \approx 1.6266,\quad
\lambda_3 \approx -2.2956,\quad
\lambda_4 \approx -0.4209.
\]
By Theorem~\ref{thm:char_polynomial}, the spectrum of 
$B_4(2,3,1,2)$ is
$\{\lambda_1,\lambda_2,\lambda_3,\lambda_4,-1,-1,-1,-1\}$.
\end{example}

Theorem~\ref{thm:char_polynomial} shows that all eigenvalues of $B_t({n_1}, \ldots, {n_t})$ different from $-1$ are precisely the eigenvalues of the matrix $M$. Hence, the squared sum of all eigenvalues other than $-1$ is given by
\[
    \sum_{\lambda \neq -1} \lambda^2 = \operatorname{tr}(M^2).
\]
We therefore focus on computing $\operatorname{tr}(M^2)$ in closed form.

\subsection{Closed form of $\operatorname{tr}(M^2)$}

From now on, we divide the discussion into two cases according to the parity of the order of $M$, that is, whether $t$ is even or odd. We write
$
    t = 2s $ or $ t = 2s+1, s \in \mathbb{N}
$
:
\begin{itemize}
    \item For $t = 2s$,
    \[
        B_{2s}(n_1, \ldots, n_s, n_{s+1}, \ldots, n_{2s})
    \]
    \item For $t = 2s+1$,
    \[
        B_{2s+1}(n_1, \ldots, n_s, n_{s+1}, \ldots, n_{2s}, n_{2s+1})
    \]
\end{itemize}


First, we write the matrix $M$ for
$B_{2s}(n_1, \ldots, n_s, n_{s+1}, \ldots, n_{2s}).$
\begin{equation*}
M=\begin{bmatrix}
n_{1}-1 & n_2 & \ldots & n_s & 0 & 0 & \ldots & 0  \\
n_1 & n_{2}-1 & \ldots & n_s & 0 & 0 & \ldots & n_{2s}  \\
\vdots & \vdots & \ddots & \vdots & \vdots & \vdots & \ddots & \vdots\\
n_1 & n_2 & \ldots & n_{s}-1 & 0 & n_{s+2} & \ldots & n_{2s} \\
0 & 0 & \ldots & 0 & n_{s+1}-1 & n_{s+2} & \ldots & n_{2s} \\
0 & 0 & \ldots & n_s & n_{s+1} & n_{s+2}-1 & \ldots & n_{2s}  \\
\vdots & \vdots & \ddots & \vdots & \vdots & \vdots & \ddots & \vdots  \\
0 & n_2 & \ldots & n_s & n_{s+1} & n_{s+2} & \ldots & n_{2s}-1 
\end{bmatrix}
\end{equation*}
Since $\operatorname{tr}(M^2) = \sum_{i=1}^{2s} \lambda_i^2$, where $\lambda_1,\dots,\lambda_{2s}$ are the eigenvalues of $M$, our next result gives a closed form for $\operatorname{tr}(M^2)$ using the structural description of $B_{2s}(n_1,\ldots,n_{2s})$.


\begin{theorem}\label{trace_even}
Let $M$ be the $2s \times 2s$ matrix as defined in Theorem~\ref{thm:char_polynomial}. Then
\begin{equation}\label{eq:13}
       \operatorname{tr}(M^2) = \sum_{i=1}^{2s} \lambda_i^2 = \left( \sum_{i=1}^{s} n_i \right)^2 + \left( \sum_{i=s+1}^{2s} n_i \right)^2 + 2 \sum_{i=2}^{s} \sum_{j=2s-i+2}^{2s} n_{i} n_{j} - 2n + 2s.
   \end{equation}
\end{theorem}

\begin{proof}
Since $\operatorname{tr}(M^2)=\sum_{i,k} m_{ik}m_{ki}$. Using the definition of $M$,
\begin{equation*}
    \operatorname{tr}(M^2)
=\sum_{i=1}^{2s}(n_i-1)^2
\;+\;\sum_{i\neq j} a_{ij}a_{ji}\,n_in_j.
\end{equation*}
Since $A(G_{2s})$ is symmetric, $a_{ij}=a_{ji}\in\{0,1\}$, so
\begin{equation}\label{eq:111}
    \operatorname{tr}(M^2)
=\sum_{i=1}^{2s}(n_i-1)^2 \;+\; 2\!\!\sum_{1\le i<j\le 2s}\!\! a_{ij}\,n_in_j.
\end{equation}

Let us now simplify the last term of \eqref{eq:111}. To do so, we partition the edge set of $G_{2s}$ according to the block structure of its adjacency matrix.
The subgraph of $G_{2s}$ induced by the first $s$ vertices is a clique thus
\begin{equation*}
    2\sum_{1\le i<j\le s} n_i n_j
    = \Bigl(\sum_{i=1}^{s} n_i\Bigr)^2 - \sum_{i=1}^{s} n_i^2.
\end{equation*}
Similarly, the subgraph of $G_{2s}$ induced by the last $s$ vertices is a clique and 
\begin{equation*}
    2\sum_{s+1\le i<j\le 2s} n_i n_j
    = \Bigl(\sum_{i=s+1}^{2s} n_i\Bigr)^2 - \sum_{i=s+1}^{2s} n_i^2.
\end{equation*}
The remaining edges in $G_{2s}$ join the first and last $s$ vertices as in Definition \ref{def_Gn}. These edges (say crossing edges) are 
\begin{equation*}
    E_{\mathrm{cross}}
    = \{(i,j): 2 \le i \le s,\;\; 2s-i+2 \le j \le 2s\},
\end{equation*}
and therefore they contribute
\begin{equation*}
    2\sum_{i=2}^{s} \;\sum_{j=\,2s-i+2}^{\,2s} n_i n_j.
\end{equation*}

Substituting these three pieces into \eqref{eq:111} gives
\begin{align*}
    \operatorname{tr}(M^2)
&=\sum_{i=1}^{2s}(n_i-1)^2 \;+\Bigl(\sum_{i=1}^{s} n_i\Bigr)^2 - \sum_{i=1}^{s} n_i^2+\Bigl(\sum_{i=s+1}^{2s} n_i\Bigr)^2 - \sum_{i=s+1}^{2s} n_i^2+ 2\sum_{i=2}^{s} \;\sum_{j=\,2s-i+2}^{\,2s} n_i n_j\\
&= \Bigl(\sum_{i=1}^{s}n_i\Bigr)^2
 + \Bigl(\sum_{i=s+1}^{2s}n_i\Bigr)^2
 + 2\sum_{i=2}^{s} \sum_{j=\,2s-i+2}^{\,2s} n_i n_j+ \sum_{i=1}^{2s}(n_i-1)^2 - \sum_{i=1}^{2s}n_i^2\\
 &= \Bigl(\sum_{i=1}^{s}n_i\Bigr)^2
 + \Bigl(\sum_{i=s+1}^{2s}n_i\Bigr)^2
 + 2\sum_{i=2}^{s} \sum_{j=\,2s-i+2}^{\,2s} n_i n_j+ \sum_{i=1}^{2s}(n_i^2-2n_i+1) - \sum_{i=1}^{2s}n_i^2
\end{align*}

and simplifying the last two terms yields $-2n+2s$. This yields the claimed formula.
\end{proof}

Now consider the odd case. For a graph of type $B_{2s+1}(n_1, \ldots, n_s, n_{s+1}, \ldots, n_{2s}, n_{2s+1})$, the matrix $M$ is written as follows,
\begin{equation*}
M=\begin{bmatrix}
n_{1}-1 & n_2 & \ldots & n_s & 0 & 0 & \ldots & 0 & n_{2s+1} \\
n_1 & n_{2}-1 & \ldots & n_s & 0 & 0 & \ldots & n_{2s} & n_{2s+1} \\
\vdots & \vdots & \ddots & \vdots & \vdots & \vdots & \ddots & \vdots & \vdots \\
n_1 & n_2 & \ldots & n_{s}-1 & 0 & n_{s+2} & \ldots & n_{2s} & n_{2s+1} \\
0 & 0 & \ldots & 0 & n_{s+1}-1 & n_{s+2} & \ldots & n_{2s} & n_{2s+1} \\
0 & 0 & \ldots & n_s & n_{s+1} & n_{s+2}-1 & \ldots & n_{2s} & n_{2s+1} \\
\vdots & \vdots & \ddots & \vdots & \vdots & \vdots & \ddots & \vdots & \vdots \\
0 & n_2 & \ldots & n_s & n_{s+1} & n_{s+2} & \ldots & n_{2s}-1 & n_{2s+1} \\
n_1 & n_2 & \ldots & n_s & n_{s+1} & n_{s+2} & \ldots & n_{2s} & n_{2s+1}-1
\end{bmatrix}
\end{equation*}

\begin{theorem}\label{trace_odd} 
Let $M$ be the $(2s+1) \times (2s+1)$ matrix as defined in Theorem~\ref{thm:char_polynomial}. Then
\begin{align*}
    \operatorname{tr}(M^2)
=\sum_{i=1}^{2s+1}\lambda_i^2=\Bigl(\sum_{i=1}^{s}n_i\Bigr)^2
 + \Bigl(\sum_{i=s+1}^{2s}n_i\Bigr)^2
 &+ 2\sum_{i=2}^{s}\sum_{j=\,2s-i+2}^{\,2s} n_i n_j \notag
 \\
&- n_{2s+1}^2
 + 2n\,n_{2s+1}
 + (2s+1) - 2n.
\end{align*}
\end{theorem}

\begin{proof}
Again we have
\begin{align}
    \operatorname{tr}(M^2)=\sum_{i,k} m_{ik}m_{ki}
&=\sum_{i=1}^{t}(n_i-1)^2 \;+\; \sum_{i\ne j} a_{ij}a_{ji}\,n_in_j \notag \\
&=\sum_{i=1}^{2s+1}(n_i-1)^2 \;+\; 2\sum_{1\le i<j\le 2s+1}\!\! a_{ij}\,n_in_j.
\label{eq:14}
\end{align}

Let us now simplify the last term of \eqref{eq:14}. We again partition the edge set of $G_{2s+1}$ according to the block structure of its adjacency matrix. 
The subgraph of $G_{2s+1}$ induced by the first $s$ vertices is a clique thus
\begin{equation*}
    2\sum_{1\le i<j\le s} n_i n_j
    = \Bigl(\sum_{i=1}^{s} n_i\Bigr)^2 - \sum_{i=1}^{s} n_i^2.
\end{equation*}

Similarly, the subgraph of $G_{2s+1}$ induced by the vertices indexed by $\{s+1,\dots,2s\}$ is a clique and 
\begin{equation*}
    2\sum_{s+1\le i<j\le 2s} n_i n_j
    = \Bigl(\sum_{i=s+1}^{2s} n_i\Bigr)^2 - \sum_{i=s+1}^{2s} n_i^2.
\end{equation*}

The edges between these two halves (the first $s$ and the next $s$ vertices) follow the same pattern as in $G_{2s+1}$. Thus the set of crossing edges is
\begin{equation*}
    E_{\mathrm{cross}}
    = \bigl\{(i,j): 2 \le i \le s,\;\; 2s - i + 2 \le j \le 2s\bigr\},
\end{equation*}
and therefore they contribute
\begin{equation*}
    2\sum_{i=2}^{s} \;\sum_{j=\,2s-i+2}^{\,2s} n_i n_j.
\end{equation*}

Finally, the vertex indexed by $2s+1$ in $G_{2s+1}$, is adjacent to every other vertex. Therefore the edges incident with part $2s+1$ contribute
\begin{equation*}
    2\sum_{i=1}^{2s} n_{2s+1} n_i
    = 2n\,n_{2s+1} - 2n_{2s+1}^2,
\end{equation*}

Substituting these four contributions into the edge–sum term of~\eqref{eq:14} and then adding the diagonal
contribution
\[
\sum_{i=1}^{2s+1} (n_i - 1)^2
= \sum_{i=1}^{2s+1} (n_i^2 - 2n_i + 1),
\]
we obtain
\begin{align*}
\operatorname{tr}(M^2)
& = \sum_{i=1}^{2s+1} (n_i^2 - 2n_i + 1)+ \Bigl(\sum_{i=1}^{s} n_i\Bigr)^2 - \sum_{i=1}^{s} n_i^2 + \Bigl(\sum_{i=s+1}^{2s} n_i\Bigr)^2 - \sum_{i=s+1}^{2s} n_i^2  \\
&\quad\quad\hspace{3cm}+ 2\sum_{i=2}^{s}\;\sum_{j=\,2s-i+2}^{\,2s} n_i n_j + 2n\,n_{2s+1} - 2n_{2s+1}^2 \\
&= \Bigl(\sum_{i=1}^{s} n_i\Bigr)^2
   + \Bigl(\sum_{i=s+1}^{2s} n_i\Bigr)^2
   + 2\sum_{i=2}^{s} \;\sum_{j=\,2s-i+2}^{\,2s} n_i n_j \\
&\quad\quad\hspace{3cm}
   + \sum_{i=1}^{2s+1} (n_i^2 - 2n_i + 1)
   - \sum_{i=1}^{2s} n_i^2
   + 2n\,n_{2s+1} - 2n_{2s+1}^2.
\end{align*}

The square–sum terms simplify to
\[
\sum_{i=1}^{2s+1} n_i^2 - \sum_{i=1}^{2s} n_i^2 - 2n_{2s+1}^2= -\,n_{2s+1}^2,
\]
and 
\[
\sum_{i=1}^{2s+1}(-2n_i + 1) = -2n + 2s+1.
\]
Substituting these back yields exactly
\[
\operatorname{tr}(M^2)
= \Bigl(\sum_{i=1}^{s} n_i\Bigr)^2
  + \Bigl(\sum_{i=s+1}^{2s} n_i\Bigr)^2
  + 2\sum_{i=2}^{s} \;\sum_{j=\,2s-i+2}^{\,2s} n_i n_j
  \;-\; n_{2s+1}^2
  \;+\; 2n\,n_{2s+1}
  \;+\; (2s+1) - 2n,
\]
This yields the claimed formula.
\end{proof}

Now, in the next subsection, we compare a structural lower bound on the clique number $\omega(G)$ with a spectral lower bound derived using the eigenvalues of $G$.

\subsection{Lower Bounds on $\omega(G)$: Structural vs.\ Spectral}

Using the structure of graph $B_{t}(n_1, \ldots, n_{t})$, we obtain the following
 bound on the clique number $\omega(G)$. For even $t=2s$, the clique number of graph
$G = B_{2s}(n_1, \ldots, n_s, n_{s+1}, \ldots n_{2s})$
 is 
\[
\omega(G)
\geq \max\left\{\, \sum_{i=1}^s n_i, \sum_{i=s+1}^{2s} n_i \,\right\}.
\]
For odd $t=2s+1$, the clique number of graph 
$G = B_{2s+1}(n_1, \ldots, n_s, n_{s+1}, \ldots, n_{2s}, n_{2s+1})$ is
\[
\omega(G)
\geq \max\left\{\, \sum_{i=1}^s n_i, \sum_{i=s+1}^{2s} n_i \,\right\} + n_{2s+1}.
\]
In particular, in both cases we have the structural bound
\[
\omega(G) \;\ge\; \frac{n}{2},
\]
where $n = \sum_{i=1}^t n_i$ is the number of vertices in $G$.

Now, we define 
\[
D := \sum_{i=1}^{t} \lambda_i^2,
\]
here $\lambda_i$ are eigenvalues of matrix $M$ corresponding to graph $B_{t}(n_1, \ldots, n_{t})$. We will prove a spectral lower bound
\[
\omega(G) \;\ge\; \frac{n}{\,n - \sqrt{D}\,}.
\]
Thus, we have two lower bounds on the clique number:
\[
\text{(A)}\quad \omega(G)\ \ge\ \frac{n}{2},
\qquad\qquad
\text{(B)}\quad \omega(G)\ \ge\ \frac{n}{\,n-\sqrt{D}\,}.
\]

A natural question is: which bound is stronger? Observation~\ref{obs:incomparable} (given below) shows that they are, in fact, incomparable.
To show the bounds are incomparable, we give examples where each bound strictly dominates the other. 




\begin{observation}\label{obs:incomparable}
For the family of graphs $B_t(n_1,\hdots, n_t)$, the bounds \textnormal{(A)} and \textnormal{(B)} are incomparable:
there exist graphs for which \textnormal{(A)} is stronger and graphs for which \textnormal{(B)} is stronger.
\end{observation}

We elaborate on this observation. To decide when (B) is no better than (A), compare the right-hand sides:
\begin{equation}\label{eq:11}
    \frac{n}{\,n-\sqrt{D}\,} \ \le\ \frac{n}{2}
\quad\Longleftrightarrow\quad
\sqrt{D}\ \le\ n-2
\quad\Longleftrightarrow\quad
D\ \le\ (n-2)^2.
\end{equation}
Substituting $s=2$ in \eqref{eq:13} and simplifying, yields
\[
D
= (n_1+n_2)^2+(n_3+n_4)^2+2\,n_2n_4-2n+4.
\]
Hence the condition $D\le(n-2)^2$ by \eqref{eq:11} is equivalent to
\begin{align}
(n_1+n_2)^2+(n_3+n_4)^2+2\,n_2n_4-2n+4&\ \le\ n^2-4n+4,\notag\\[2pt]
-2(n_1+n_2)(n_3+n_4)+2\,n_2n_4+2n&\ \le\ 0,\notag\\[2pt]
n - n_1n_3 - n_1n_4 - n_2n_3&\ \le\ 0. \label{eq:12}
\end{align}
Therefore:
\begin{itemize}
\item If \eqref{eq:12} holds, then 
$\dfrac{n}{n-\sqrt{D}}\le\dfrac{n}{2}$ and bound \textnormal{(A)} is at least as strong as \textnormal{(B)}.
\item If \eqref{eq:12} fails, then 
$\dfrac{n}{n-\sqrt{D}}>\dfrac{n}{2}$ and bound \textnormal{(B)} is strictly stronger than \textnormal{(A)}.
\end{itemize}

It remains to give concrete valid choices of $(n_1,n_2,n_3,n_4)$
exhibiting both regimes. Take $n=8$ in both examples:
\begin{itemize}
\item If $(n_1,n_2,n_3,n_4)=(1,2,2,3)$, then
\[
n - n_1n_3 - n_1n_4 - n_2n_3
= 8 - (1\cdot2) - (1\cdot3) - (2\cdot2) = -1 \le 0,
\]
so \eqref{eq:12} holds and (A) is stronger.
\item If $(n_1,n_2,n_3,n_4)=(1,2,1,4)$, then
\[
n - n_1n_3 - n_1n_4 - n_2n_3
= 8 - (1\cdot1) - (1\cdot4) - (2\cdot1) = 1 > 0,
\]
so \eqref{eq:12} fails and (B) is stronger.
\end{itemize}
These two instances show that the bounds (A) and (B) are incomparable. 


\section{Main Result}\label{sec:main_result}
In this section we consider the families of graphs $B_3(n_1,n_2,n_3)$ and
$B_4(n_1,n_2,n_3,n_4)$.
For each such graph we take the associated $t \times t$ matrix $M$ and denote its eigenvalues by
$\lambda_1,\dots,\lambda_t$.
These are precisely the adjacency eigenvalues of $B_t(n_1,\dots,n_t)$ other than $-1$.

Our goal is to establish the Wilf-type inequality
\[
\sqrt{\sum_{i=1}^{t}\lambda_i^2}\;\le\; n\!\left(1-\frac{1}{\omega}\right),
\]
for $B_t(n_1,\dots,n_t)$.
First, we prove the inequality for $B_3(n_1,n_2,n_3)$ in Theorem \ref{thm:omega-odd-3}. 



\begin{theorem}\label{thm:omega-odd-3}
Let $G=B_{3}(n_1, n_2, n_3)$ be a graph of order $n$ with positive integers $n_1,n_2,n_3$. Let $\lambda_1,\lambda_2,\lambda_3$ be the adjacency eigenvalues of $G$ different from $-1$. If $(n_1,n_2)\neq (1,1)$, then
\begin{equation}\label{eq:omega-odd-3}
\sqrt{\lambda_1^2+\lambda_2^2+\lambda_3^2}\ \le\ n \left(1-\frac{1}{\omega}\right).
\end{equation}
In particular, for $n_1=n_2=1$ and $n_3\ge 1$, the inequality \eqref{eq:omega-odd-3} does not hold.
\end{theorem}

\begin{proof}
Define
\begin{equation*}\label{eq:def-pqn-omega}
n:=n_1+n_2+n_3,\qquad \omega:=\max\{n_1+n_3,n_2+n_3\}.
\end{equation*}
We consider three cases: $n_1>n_2$, $n_1<n_2$, and $n_1=n_2$.\\

\noindent\textbf{Case 1: $n_1>n_2$.}\\
By Theorem~\ref{trace_odd} with $s=1$,
\begin{equation}\label{eq:case1-trace}
\lambda_1^2+\lambda_2^2+\lambda_3^2
= n_1^2+n_2^2 - n_3^2 + 2n\,n_3 + 3 - 2n,
\end{equation}
Since $\omega=n_1+n_3$ in this case, squaring
\eqref{eq:omega-odd-3} and using $\lambda_1^2+\lambda_2^2+\lambda_3^2\ge 0$,
we obtain the equivalent inequality
\begin{equation}\label{eq:case1-goal-squared}
\lambda_1^2+\lambda_2^2+\lambda_3^2
\ \le\ 
n^2\left(1-\frac{1}{n_1+n_3}\right)^2 .
\end{equation}
Substituting $n=n_1+n_2+n_3$ into \eqref{eq:case1-trace} gives
\begin{equation}\label{eq:case1-sumlambda-expanded}
\lambda_1^2+\lambda_2^2+\lambda_3^2
=
n_1^2+n_2^2-n_3^2+2(n_1+n_2+n_3)n_3+3-2(n_1+n_2+n_3).
\end{equation}

Define
\begin{equation*}\label{eq:case1-Delta-def}
\Delta(n_1,n_2,n_3)
:=
n^2\left(1-\frac{1}{n_1+n_3}\right)^2
-\bigl(\lambda_1^2+\lambda_2^2+\lambda_3^2\bigr).
\end{equation*}
Then \eqref{eq:case1-goal-squared} is equivalent to
\begin{equation}\label{eq:case1-Delta-nonneg}
\Delta(n_1,n_2,n_3)\ge 0.
\end{equation}
Clearing denominators using \eqref{eq:case1-sumlambda-expanded} yields
\begin{equation}\label{eq:case1-Delta-D}
\Delta(n_1,n_2,n_3)=\frac{D(n_1,n_2,n_3)}{(n_1+n_3)^2},
\end{equation}
where
\begin{equation}\label{eq:case1-D-def}
\begin{aligned}
D(n_1,n_2,n_3)
=\ &(n_1+n_3-1)^2(n_1+n_2+n_3)^2\\
&-(n_1+n_3)^2\Bigl(n_1^2+n_2^2-n_3^2+2(n_1+n_2+n_3)n_3+3-2(n_1+n_2+n_3)\Bigr).
\end{aligned}
\end{equation}
Since $(n_1+n_3)^2>0$, by \eqref{eq:case1-Delta-D} the condition
\eqref{eq:case1-Delta-nonneg} is equivalent to
\begin{equation}\label{eq:case1-D-nonneg}
D(n_1,n_2,n_3)\ge 0.
\end{equation}

Now fix $n_1$ and $n_3$ and write $D(n_2):=D(n_1,n_2,n_3)$ for $n_2\in[1,n_1-1]$.
Set
\begin{equation*}\label{eq:case1-a}
a:=n_1+n_3.
\end{equation*}
A direct expansion of \eqref{eq:case1-D-def} shows that $D(n_2)$ is a quadratic
polynomial in the single variable $n_2$:
\begin{equation*}\label{eq:case1-D-quadratic}
D(n_2)=(1-2a)\,n_2^2+\beta(n_1,n_3)\,n_2-2a^2,
\end{equation*}
where
\begin{equation*}\label{eq:case1-beta}
\beta(n_1,n_3)
:=2n_1^3+4n_1^2n_3-2n_1^2+2n_1n_3^2-4n_1n_3+2n_1-2n_3^2+2n_3 .
\end{equation*}
Since $n_1,n_3\ge 1$, we have $a\ge 2$, hence
\begin{equation*}\label{eq:case1-concave}
1-2a<0,
\end{equation*}
so $D(n_2)$ is concave down on $\mathbb{R}$. Therefore, on the interval $[1,n_1-1]$ its
minimum is attained at an endpoint:
\begin{equation}\label{eq:case1-endpoint-reduction}
\min_{\,n_2\in[1,n_1-1]} D(n_2)=\min\{D(1),D(n_1-1)\}.
\end{equation}
It thus suffices to prove
\begin{equation*}\label{eq:case1-endpoints-suffice}
D(1)\ge 0
\qquad\text{and}\qquad
D(n_1-1)\ge 0.
\end{equation*}

\smallskip
Now substituting $n_2=1$ into \eqref{eq:case1-D-def} and simplifying yields
\begin{equation*}\label{eq:case1-D1}
D(1)=D(n_1,1,n_3)=2(n_1-2)(n_1+n_3)^2+1=2(n_1-2)a^2+1.
\end{equation*}
Since $n_1>n_2\ge 1$, we have $n_1\ge 2$. If $n_1=2$, then $D(1)=1$.
If $n_1\ge 3$, then $n_1-2\ge 1$ and $a^2>0$, hence $D(1)\ge 1$.
Thus
\begin{equation}\label{eq:case1-D1-nonneg}
D(1)\ge 1.
\end{equation}

\smallskip
Substituting $n_2=n_1-1$ into \eqref{eq:case1-D-def} and simplifying gives
\begin{equation}\label{eq:case1-Dn1m1}
D(n_1-1)=D(n_1,n_1-1,n_3)
=(n_1-1)^2+2(n_1-2)a\bigl(n_1a-(n_1-1)\bigr).
\end{equation}
If $n_1=2$, then $n_2=1$ and $D(n_1-1)=D(1)=1$.
Assume $n_1\ge 3$. Then $n_1-2>0$ and $a=n_1+n_3\ge n_1+1\ge 4$, and
\begin{equation*}\label{eq:case1-bracket-pos}
n_1a-(n_1-1)=n_1(a-1)+1\ge 1.
\end{equation*}
Hence the second term in \eqref{eq:case1-Dn1m1} is nonnegative, so
\begin{equation}\label{eq:case1-Dn1m1-nonneg}
D(n_1-1)\ge (n_1-1)^2\ge 1.
\end{equation}

Combining \eqref{eq:case1-endpoint-reduction}, \eqref{eq:case1-D1-nonneg},
and \eqref{eq:case1-Dn1m1-nonneg}, we obtain $D(n_2)\ge 0$ for all
$n_2\in[1,n_1-1]$. Therefore \eqref{eq:case1-D-nonneg} holds, which implies
\eqref{eq:case1-Delta-nonneg} by \eqref{eq:case1-Delta-D}, hence
\eqref{eq:case1-goal-squared} and thus \eqref{eq:omega-odd-3}.
This completes the proof of Case~1.

\smallskip
\noindent\textbf{Case 2: $n_1<n_2$.}\\
This case is symmetric to Case~1 after interchanging the roles of $n_1$ and $n_2$, so the same argument proves \eqref{eq:omega-odd-3}.

\smallskip
\noindent\textbf{Case 3: $n_1=n_2$.}\\
Set $t:=n_1=n_2$. 
Recall $n=n_1+n_2+n_3=2t+n_3$ and
\begin{equation*}\label{eq:case2-omega}
\omega=\max\{n_1+n_3,n_2+n_3\}=t+n_3.
\end{equation*}
By Theorem~\ref{trace_odd} with $s=1$,
\begin{equation}\label{eq:case2-trace}
\lambda_1^2+\lambda_2^2+\lambda_3^2
= n_1^2+n_2^2-n_3^2+2nn_3+3-2n.
\end{equation}
Substituting $n_1=n_2=t$ and $n=2t+n_3$ into \eqref{eq:case2-trace} yields
\begin{equation}\label{eq:case2-sumlambda}
\lambda_1^2+\lambda_2^2+\lambda_3^2
= n_3^2+4tn_3-2n_3+2t^2-4t+3.
\end{equation}

Since $\omega\geq t+n_3$, squaring \eqref{eq:omega-odd-3} and using
$\lambda_1^2+\lambda_2^2+\lambda_3^2\ge 0$ gives the equivalent inequality
\begin{equation}\label{eq:case2-goal-squared}
\lambda_1^2+\lambda_2^2+\lambda_3^2
\ \le\
n^2\left(1-\frac{1}{t+n_3}\right)^2.
\end{equation}
Define
\begin{equation*}\label{eq:case2-Delta}
\Delta(t,n_3)
:=
n^2\left(1-\frac{1}{t+n_3}\right)^2
-\bigl(\lambda_1^2+\lambda_2^2+\lambda_3^2\bigr),
\qquad n=2t+n_3.
\end{equation*}
Then \eqref{eq:case2-goal-squared} is equivalent to
\begin{equation}\label{eq:case2-Delta-nonneg}
\Delta(t,n_3)\ge 0.
\end{equation}
Clearing denominators (using \eqref{eq:case2-sumlambda}) gives
\begin{equation}\label{eq:case2-Delta-D}
\Delta(t,n_3)=\frac{D(t,n_3)}{(t+n_3)^2},
\end{equation}
where
\begin{equation}\label{eq:case2-D-def}
D(t,n_3)
:=(t+n_3-1)^2(2t+n_3)^2-(t+n_3)^2\bigl(n_3^2+4tn_3-2n_3+2t^2-4t+3\bigr).
\end{equation}
Since $(t+n_3)^2>0$, by \eqref{eq:case2-Delta-D} it suffices to prove
\begin{equation}\label{eq:case2-D-nonneg}
D(t,n_3)\ge 0 \qquad (t\ge 2,\ n_3\ge 1).
\end{equation}

A direct expansion of \eqref{eq:case2-D-def} shows that $D(t,n_3)$
is a quadratic polynomial in the single variable $n_3$:
\begin{equation*}\label{eq:case2-D-quadratic}
D(t,n_3)=A(t)\,n_3^2+B(t)\,n_3+C(t),
\end{equation*}
where
\begin{equation*}\label{eq:case2-ABC}
A(t)=2(t^2-t-1),\qquad
B(t)=2t(2t^2-3t-1),\qquad
C(t)=t^2(2t^2-4t+1).
\end{equation*}
For $t\ge 2$, we have $A(t)>0$ and $B(t)>0$. Hence for every $n_3\ge 1$,
\begin{equation*}\label{eq:case2-monotone}
D(t,n_3)-D(t,1)=A(t)(n_3^2-1)+B(t)(n_3-1)\ \ge\ 0.
\end{equation*}
So the minimum over $n_3\ge 1$ is attained at $n_3=1$:
\begin{equation}\label{eq:case2-reduction-n3}
D(t,n_3)\ge D(t,1)\qquad (t\ge 2,\ n_3\ge 1).
\end{equation}
Substituting $n_3=1$ into \eqref{eq:case2-D-def} gives
\begin{equation*}\label{eq:case2-D-at-1}
D(t,1)=2t^4-3t^2-4t-2.
\end{equation*}
Using $2(t-1)^2\ge 0$, we have $-4t\ge -2t^2-2$, and therefore
\begin{equation*}\label{eq:case2-lowerbound}
D(t,1)=2t^4-3t^2-4t-2\ \ge\ 2t^4-5t^2-4.
\end{equation*}
The polynomial $\phi(t):=2t^4-5t^2-4$ is increasing for $t\ge 2$ since
\begin{equation*}\label{eq:case2-phi-derivative}
\phi'(t)=8t^3-10t=2t(4t^2-5)>0\qquad (t\ge 2),
\end{equation*}
and $\phi(2)=8>0$. Hence
\begin{equation}\label{eq:case2-D1-positive}
D(t,1)\ge 0 \qquad (t\ge 2).
\end{equation}
Combining \eqref{eq:case2-reduction-n3} and \eqref{eq:case2-D1-positive} yields
$D(t,n_3)\ge 0$ for all $t\ge 2$ and $n_3\ge 1$. Consequently,
\eqref{eq:case2-D-nonneg} holds, which implies \eqref{eq:case2-Delta-nonneg} by
\eqref{eq:case2-Delta-D}, and therefore \eqref{eq:case2-goal-squared} and
\eqref{eq:omega-odd-3}. This completes the proof of Case~3 when $t\ge 2$.

Now consider $t=1$. Then $n_1=n_2=1$ and $n_3\ge 1$, and hence
\begin{equation}\label{eq:t1-n-omega}
n=n_1+n_2+n_3=n_3+2,
\qquad
\omega=\max\{1+n_3,1+n_3\}=n_3+1.
\end{equation}
Substituting $t=1$ into \eqref{eq:case2-sumlambda} gives
\begin{equation*}\label{eq:t1-sumsq}
\lambda_1^2+\lambda_2^2+\lambda_3^2
=n_3^2+4n_3-2n_3+2-4+3
=(n_3+1)^2,
\end{equation*}
and hence
\begin{equation}\label{eq:t1-root}
\sqrt{\lambda_1^2+\lambda_2^2+\lambda_3^2}=n_3+1.
\end{equation}
We can rewrite \eqref{eq:omega-odd-3} as
\begin{equation}\label{eq:omega-odd-31}
\frac{n}{\,n-\sqrt{\lambda_1^2+\lambda_2^2+\lambda_3^2}\,}\ \le\ \omega.
\end{equation}
Substituting \eqref{eq:t1-n-omega} and \eqref{eq:t1-root} into \eqref{eq:omega-odd-31} yields
\begin{equation*}\label{eq:t1-contradiction-chain}
\frac{n_3+2}{(n_3+2)-(n_3+1)}\le n_3+1
\quad\Longrightarrow\quad
n_3+2\le n_3+1,
\end{equation*}
which is a contradiction. Therefore \eqref{eq:omega-odd-3} fails for $B_3(1,1,n_3)$.


This completes the proof.
\end{proof}

\begin{example}
Consider the graph $B_3(2,3,2)$ of order $n=7$ and the associated matrix $M$ of order $t=3$:
\begin{equation*}
M=
\begin{bmatrix}
1&0&2\\
0&2&2\\
2&3&2
\end{bmatrix}.
\end{equation*}
Let $\lambda_1,\lambda_2,\lambda_3$ denote the eigenvalues of $M$. A direct computation gives
\begin{equation}\label{eq:exB3-232-eigsM11}
\lambda_1 \approx 4.52,\quad
\lambda_2 \approx 1.39,\quad
\lambda_3 \approx -1.91.
\end{equation}
Hence, by Theorem~\ref{thm:char_polynomial}, the eigenvalues of $B_3(2,3,2)$ are $-1$ with multiplicity $n-t=4$ together with \eqref{eq:exB3-232-eigsM11}, i.e.,
\begin{equation*}
\operatorname{spec}\bigl(B_3(2,3,2)\bigr)
=\{4.52,\ 1.39,\ -1.91,\ -1,\ -1,\ -1,\ -1\}.
\end{equation*}
Therefore
\begin{equation*}
\sum_{i=1}^{3}\lambda_i^2
\approx (4.52)^2+(1.39)^2+(1.91)^2
\approx 26,
\qquad\text{so}\qquad
\sqrt{\sum_{i=1}^{3}\lambda_i^2}\approx 5.10.
\end{equation*}
Since $\omega=\max\{n_1+n_3,\;n_2+n_3\}=\max\{4,5\}=5$, we have
\begin{equation*}
n\!\left(1-\frac{1}{\omega}\right)=7\!\left(1-\frac{1}{5}\right)=\frac{28}{5}=5.60,
\end{equation*}
and hence the inequality \eqref{eq:omega-odd-3}
holds for $B_3(2,3,2)$.
\end{example}

\begin{remark}\label{rem:B3-11n3}
We now verify that the Elphick--Wocjan conjecture holds for $B_3(1,1,n_3)$. Note that $B_3(1,1,n_3)$ is the complete multipartite graph $K_{2,1,1,\ldots,1}$ and consequently it has exactly one positive eigenvalue \cite{Smith1970}. Then $\sqrt{s^+}=\lambda_1$. Thus for this
family the Elphick--Wocjan conjecture is equivalent to the Wilf bound, which always holds.
\end{remark}

The Elphick--Wocjan conjecture follows  from Theorem~\ref{thm:omega-odd-3} because
$$\sum_{\lambda_i\neq -1}\lambda_i^{2}\;\ge\;\sum_{\lambda_i>0}\lambda_i^{2}\;=\;s^{+},$$
and the $n_1=n_2=1$ case discussed in Remark~\ref{rem:B3-11n3}.

\begin{corollary}\label{co:11}
The Elphick-Wocjan conjecture holds for \(B_3(n_1, n_2, n_3) \) for any positive integers \( n_1, n_2, n_3 \).  
\end{corollary}


\bigskip
Now we investigate the infinite family of graphs \(B_4(n_1, n_2, n_3, n_4) \). First note for \( B_4(1,1,1,1)=P_4 \)  that 
\[
\lambda_1=\frac{1+\sqrt5}{2},\quad
\lambda_2=\frac{\sqrt5-1}{2},\quad
\lambda_3=-\frac{\sqrt5-1}{2},\quad
\lambda_4=-\frac{1+\sqrt5}{2}.
\]
Therefore
\[
\sum_{i=1}^{4}\lambda_i^2
=2\left(\frac{1+\sqrt5}{2}\right)^2+2\left(\frac{\sqrt5-1}{2}\right)^2
=6,
\qquad\text{so}\qquad
\sqrt{\sum_{i=1}^{4}\lambda_i^2}=\sqrt6>2.
\]
Since $\omega(P_4)=2$, we have $n\!\left(1-\frac1\omega\right)=4\left(1-\frac12\right)=2$ and thus the inequality
\[
\sqrt{\sum_{i=1}^{4}\lambda_i^2}\ \le\ n\!\left(1-\frac{1}{\omega}\right)
\]
fails for \( B_4(1,1,1,1)=P_4 \). 
Now to verify the Elphick--Wocjan conjecture, $s^+=\lambda_1^2+\lambda_2^2=3$ and 
\[
\sqrt{s^+}=\sqrt3<2=n\!\left(1-\frac12\right).
\]
Thus the Elphick--Wocjan conjecture holds for $B_4(1,1,1,1)$.



\begin{theorem}
Let \( G = B_4(n_1, n_2, n_3, n_4) \) be a graph of order $n$ with positive integers \( n_1, n_2, n_3, n_4 \). Let $\lambda_1,\lambda_2,\lambda_3,\lambda_4$ be  the adjacency eigenvalues of $G$ different from $-1$. 
Then
\begin{equation}\label{eq:omega-bound}
\sqrt{\sum_{i=1}^{4}\lambda_i^2}\;\le\; n\!\left(1-\frac{1}{\omega}\right),
\end{equation}
except when
\( (n_1,n_2,n_3,n_4)=(1,1,1,1)
\).
\end{theorem}

\begin{proof}
By Theorem~\ref{thm:char_polynomial}, the matrix \( M \) has eigenvalues 
\( \lambda_1, \lambda_2, \lambda_3, \lambda_4 \). In particular,
\[
    \operatorname{tr}(M^2) = \sum_{i=1}^{4} \lambda_i^2.
\] Hence the Wilf-type bound \eqref{eq:omega-bound} is equivalent to
\begin{equation}\label{eq:trM2-omega}
    \operatorname{tr}(M^2)\;\le\; n^2\!\left(1-\frac{1}{\omega}\right)^2.
\end{equation}

We set
\[
    p := n_1 + n_2,\quad q := n_3 + n_4,\quad  \omega \geq \max\{p,q\},  \quad n= p + q=\sum_{i=1}^{4}n_i.
\]
To establish \eqref{eq:trM2-omega}, we consider the following three cases:
\[
    \text{(1) } p > q,\qquad \text{(2) } p = q,\qquad \text{(3) } p < q.
\]

\medskip
\noindent\textbf{Case 1: $p>q$.}\\
In this case $\omega\geq\max\{p,q\}=p$ and $s=2$. 
Put $s=2$ in \eqref{eq:13} to obtain \eqref{eq:trM2-omega} with $\omega \ge p$. 
It therefore suffices to prove:
\begin{equation}\label{eq:even-general}
p^2 
+ q^2 
+ 2  n_2n_4 - 2n + 4 
\;\leq\; n^2\!\left(1 - \frac{1}{p}\right)^2.
\end{equation}
Expanding the right-hand side,
\[
n^2\!\left(1 - \frac{1}{p}\right)^2
= n^2 - \frac{2n^2}{p} + \frac{n^2}{p^2}.
\]
Since
\[
n^2= (p+q)^2
= p^2 
  + q^2 
  + 2 p q,
\]
inequality \eqref{eq:even-general} can be rewritten as
\begin{align*}
p^2 
+ q^2 
+ 2  n_2n_4 - 2n + 4  \leq 
 p^2 
  + q^2 
  + 2 p q  - \frac{2n^2}{p} + \frac{n^2}{p^2}.
\end{align*}

Cancelling the common terms from both sides yields
\begin{equation}\label{eq:gen-s-ineq}
2n_2n_4 - 2n + 4 \;\le\; 
  2pq - \frac{2n^2}{p} + \frac{n^2}{p^2}.
\end{equation}
Since $n = p+q$, we rewrite \eqref{eq:gen-s-ineq} as
\begin{equation}\label{eq4.2}
2n_2n_4 - 2(p+q) + 4 \;\le\;
  2pq - 2(p+q)\left(1 + \frac{q}{p}\right)
  + \left(1 + \frac{q}{p}\right)^2.
\end{equation}
Rewriting the above inequality, we get the following

\begin{equation}
\begin{aligned}
    0
    &\leq 
    2pq - 2n_2 n_4 - 2(p+q) \frac{q}{p} + \left(1+ \frac{q}{p}\right)^2-4\\
    &= 2pq - 2n_2n_4 - 2q - 2 \frac{q^2}{p}  +  \frac{q^2}{p^2} + 2 \frac{q}{p}-3,
\end{aligned}
\end{equation}


For the worst case (to minimize the right-hand side) we make the negative term
$-2n_2n_4$ as small as possible. Since
\[
n_2 \le p-1,\qquad n_4 \le q-1,
\]
we have
\[
n_2n_4 \le (p-1)(q-1)
\quad\Longrightarrow\quad
-2n_2n_4 \;\ge\; -2(p-1)(q-1).
\]
Therefore
\begin{align*}
    &2pq - 2n_2n_4 - 2q - 2 \frac{q^2}{p}  +  \frac{q^2}{p^2} + 2 \frac{q}{p}-3\\
    &\hspace{3cm}\geq2pq - 2(p-1)(q-1)  - 2q - 2 \frac{q^2}{p}  +  \frac{q^2}{p^2} + 2 \frac{q}{p}-3\\
    &\hspace{3cm}=2pq - 2(pq-p-q+1) - 2q - 2 \frac{q^2}{p}  +  \frac{q^2}{p^2} + 2 \frac{q}{p}-3\\
    &\hspace{3cm}=2p + 2q - 2 - 2q - 2 \frac{q^2}{p}  +  \frac{q^2}{p^2} + 2 \frac{q}{p}-3\\
    &\hspace{3cm}=2p  - 2 \frac{q^2}{p}  +  \frac{q^2}{p^2} + 2 \frac{q}{p}-5.
\end{align*}


Now to complete Case~1 it suffices to show that
\begin{equation}\label{eq:Hminus4}
    2p - 2\frac{q^2}{p} + \frac{q^2}{p^2} + 2\frac{q}{p} - 5 \;\ge\; 0.
\end{equation}
for all integers $p\ge3$, $q\ge2$ with $q \le p-1$ and $p>q$.  
Since $p>0$, multiplying \eqref{eq:Hminus4} by $p^2$ preserves the inequality and gives
\[
    p^2\!\left(
      2p - 2\frac{q^2}{p} + \frac{q^2}{p^2} + 2\frac{q}{p} - 5
    \right)
    = 2p^3 - 5p^2 + 2pq - 2pq^2 + q^2.
\]
Define, for a fixed $p\geq 3$,
\[
    f_p(q)
    := 2p^3 - 5p^2 + 2pq - 2pq^2 + q^2
     = (1-2p)q^2 + 2pq + (2p^3-5p^2).
\]
Then \eqref{eq:Hminus4} is equivalent to $f_p(q)\ge 0$ for all integers
$2\le q\le p-1$. Note that $f_p(q)$ is a quadratic in $q$ with leading coefficient
$1-2p<0$ (since $p\ge3$), so $f$ opens downward. Then $f(q)$ achieve its maximum at
\[
    q_{\mathrm{v}}
    = -\frac{2p}{2(1-2p)}
    = \frac{p}{2p-1},
\]
and for $p\ge3$ we have $q_{\mathrm{v}} \le \tfrac{3}{5} < 1$.
Hence $f$ is strictly decreasing on $[2,\infty)$.

In our range $2 \le q \le p-1$, the minimum of $f$ is therefore attained at
$q=p-1$. A direct computation yields
\begin{align*}
    f(p-1)
    &=(1-2p)(p-1)^2 + 2p(p-1)+(2p^3 - 5p^2)= 2p^2 - 6p + 1.
\end{align*}
The polynomial $g(p) := 2p^2 - 6p + 1$ is upward-opening, with discriminant $\Delta = 28,$
so its roots are
\[
    p = \frac{6 \pm \sqrt{28}}{4}
      = \frac{3 \pm \sqrt{7}}{2}.
\]
Thus $g(p) \ge 0$ for all integers $p \ge 3$. Consequently
$f(q) \ge 0$
for all integers $p \ge 3$ and $2 \le q \le p-1$ and hence \eqref{eq:Hminus4} holds. 


\medskip
\noindent\textbf{Case 2: $p=q$.}\\
Here $\omega\geq \max\{p,q\}=p=q$, so $n=2p$ and $q/p=1$. Returning to \eqref{eq4.2} and
substituting $q/p=1$ and $n=2p$, the right-hand side becomes
\[
2pq - 2n\!\left(1+\frac{q}{p}\right) + \left(1+\frac{q}{p}\right)^2
= 2p^2 - 2n\cdot 2 + 2^2
= 2p^2 - 4n + 4
= 2p^2 - 8p + 4.
\]
Thus \eqref{eq4.2} is equivalent to
\begin{align*}
2 n_2n_4 - 2n + 4 \;\le\; 2p^2 - 4n + 4
&\;\Longleftrightarrow\; 2 n_2n_4 \;\le\; 2p^2 - 2n\\
&\;\Longleftrightarrow\; n_2n_4 \;\le\; p^2 - n
 = p(p-2).
\end{align*}

On the other hand, using $n_1+n_2=p$ and $n_3+n_4=q=p$ with $n_i\ge1$ we have
$1\le n_2\le p-1$ and $1\le n_4\le p-1$, so
\[
n_2n_4 \;\le\; (p-1)^2 .
\]
Comparing, we see that the inequality
\[
n_2n_4 \;\le\; p(p-2)
\]
holds  
whenever $n_2$ or $n_4$ is at most $p-2$. The only way this inequality 
fails is when 
 $n_2=n_4=p-1$ (consequently $n_1=p-n_2=1$ and $n_3=p-n_4=1$), i.e., when 
\[
(n_1,n_2,n_3,n_4)=(1,p-1,\,1,p-1).
\]
Now we show that \eqref{eq:omega-bound} holds for all \(p\ge3\)
and fails only for \(p=2\), i.e., only for \((n_1,n_2,n_3,n_4)=(1,1,1,1)\).

In this case, the clique number 
\[
\omega\geq 2p-2.
\] 
A direct substitution $n_1=n_3=1 $ and $n_2=n_4=p-1$ into
\eqref{eq:13} gives 
\[
\operatorname{tr}(M^2)=4p^2-8p+6.
\]
With $n=2p$ and $\omega\ge 2p-2$, we have
\[
n^2\!\left(1-\frac{1}{\omega}\right)^2
\ge (2p)^2\left(1-\frac{1}{2p-2}\right)^2.
\]
Expanding the left-hand side gives
\[
(2p)^2\left(1-\frac{1}{2p-2}\right)^2
=4p^2\cdot \frac{(2p-3)^2}{4(p-1)^2}
=\frac{p^2(2p-3)^2}{(p-1)^2}
=\frac{4p^4-12p^3+9p^2}{(p-1)^2}.
\]
Therefore, to prove \eqref{eq:omega-bound} it suffices to show
\[
4p^2-8p+6 \le \frac{4p^4-12p^3+9p^2}{(p-1)^2},
\]
equivalently it reduce to
\[
0\le (4p^4-12p^3+9p^2)-(4p^2-8p+6)(p-1)^2 
=4p^3-17p^2+20p-6.
\]
Set
\[
f(p):=4p^3-17p^2+20p-6.
\]
Then $f(2)=-2<0$, so the inequality fails at $p=2$, i.e., for $(n_1,n_2,n_3,n_4)=(1,1,1,1)$.
Moreover,
\[
f'(p)=12p^2-34p+20,
\]
whose real roots are \(p=\frac{10}{12}<2\) and \(p=\frac{24}{12}=2\), so
\(f'(p)>0\) for all \(p>2\). Hence \(f(p)\) is strictly increasing on
\((2,\infty)\) and since \(f(3)=9>0\), we obtain \(f(p)>0\) for every
integer \(p\ge3\). Therefore \eqref{eq:omega-bound} holds for all \(p\ge3\)
and fails only at \(p=2\), that is, only for \((n_1,n_2,n_3,n_4)=(1,1,1,1)\). This completes the proof of Case 2. 

\medskip
\noindent\textbf{Case 3: $p<q$.}\\
This case is symmetric to Case~1 after interchanging the roles of $p$ and $q$, so the same argument proves \eqref{eq:omega-bound}. This completes the proof.
\end{proof}

Since $\sum_{\lambda_i \neq -1} \lambda_i^{2}\geq\sum_{\lambda_i>0} \lambda_i^{2}$, we have the following corollary.
\begin{corollary}
The Elphick-Wocjan conjecture holds for \(B_4(n_1, n_2, n_3, n_4) \) for any positive integers \( n_1, n_2, n_3, n_4 \).
\end{corollary}

\begin{example}
Consider the graph $B_4(2,3,1,2)$ in Example~\ref{example1}. The spectrum of $B_4(2,3,1,2)$ is
\[
\{\lambda_1,\lambda_2,\lambda_3,\lambda_4,-1,-1,-1,-1\},
\]
where
\[
\lambda_1 \approx 5.09,\quad
\lambda_2 \approx 1.63,\quad
\lambda_3 \approx -2.30,\quad
\lambda_4 \approx -0.42.
\]

Therefore
\[
\sum_{i=1}^{4}\lambda_i^2
\approx (5.09)^2+(1.63)^2+(2.30)^2+(0.42)^2
\approx 34,
\qquad\text{so}\qquad
\sqrt{\sum_{i=1}^{4}\lambda_i^2}\approx 5.83.
\]
Here $n=8$ and $\omega\geq\max\{n_1+n_2,\ n_3+n_4\}=\max\{5,3\}=5$, so
\[
n\!\left(1-\frac{1}{\omega}\right)
\geq8\!\left(1-\frac{1}{5}\right)
=\frac{32}{5}=6.40.
\]
Hence
\[
\sqrt{\sum_{i=1}^{4}\lambda_i^2}\approx 5.83 \le 6.40 \leq n\!\left(1-\frac{1}{\omega}\right).
\]
\end{example}

A natural next class of graphs could be considered is $B_5(n_1,n_2,n_3,n_4,n_5)$. But we show in the following example that there are graphs other than $B_5(1,1,1,1,1)$ for which the Elphick-Wocjan conjecture holds but our strong spectral bound for the clique number is not valid.

\begin{example} Consider the graph $B_5(1,2,1,1,1)$ of order $n=6$ 
and the associated matrix $M$ of order $t=5$:
\begin{equation*}
M=
\begin{bmatrix}
0&2&0&0&1\\
1&1&0&1&1\\
0&0&0&1&1\\
0&2&1&0&1\\
1&2&1&1&0
\end{bmatrix}.
\end{equation*}
Let $\lambda_1,\dots,\lambda_5$ denote the eigenvalues of $M$. A direct computation gives
\begin{equation}\label{eq:exB5b-eigsM}
\lambda_1\approx 3.86,\quad
\lambda_2\approx 0.78,\quad
\lambda_3\approx -0.38,\quad
\lambda_4\approx -1.48,\quad
\lambda_5\approx -1.78.
\end{equation}
Hence, by Theorem~\ref{thm:char_polynomial}, the eigenvalues of $B_5(1,2,1,1,1)$ are $-1$ with multiplicity $n-t=1$ together with \eqref{eq:exB5b-eigsM}, i.e.,
\begin{equation*}
\operatorname{spec}\bigl(B_5(1,2,1,1,1)\bigr)
=\{3.86,\ 0.78,\ -0.38,\ -1,\ -1.48,\ -1.78\}.
\end{equation*}

The clique number of $B_5(1,2,1,1,1)$ is $\omega=4$. Therefore
\begin{equation*}
n\!\left(1-\frac{1}{\omega}\right)
=6\!\left(1-\frac{1}{4}\right)
=\frac{9}{2}.
\end{equation*}



For Elphick--Wocjan conjecture,
the positive eigenvalues are $\lambda_1$ and $\lambda_2$, hence
\begin{equation*}
s^+=\lambda_1^2+\lambda_2^2\approx 15.50
\end{equation*}
and thus
\begin{equation*}
\sqrt{s^+}\approx 3.94\ \le\ \frac{9}{2}=n\!\left(1-\frac{1}{\omega}\right).
\end{equation*}

For our bound,
since $\sum_{i=1}^{5}\lambda_i^2=\operatorname{tr}(M^2)$ and $\operatorname{tr}(M^2)=21$, we obtain
\begin{equation*}
\sum_{\lambda\in\operatorname{spec}(M)} \lambda^2=\operatorname{tr}(M^2)=21.
\end{equation*}
Consequently,
\begin{equation*}
\sqrt{\sum_{\lambda\in\operatorname{spec}(M)} \lambda^2}
=\sqrt{21}\approx 4.58 \not\le \frac{9}{2}=n\!\left(1-\frac{1}{\omega}\right).
\end{equation*}


In this example, the Elphick-Wocjan conjecture holds. 
However, our stronger bound fails for $B_5(1,2,1,1,1)$. 
This suggests that, for $t=5$, there may be additional parameter choices for which our bound fails.

\end{example}

\subsection*{AI Disclosure}

The authors declare that no AI tools were used to conduct this research.







\bibliographystyle{splncs04}
\bibliography{ref}

@article{wilf,
  title={Spectral bounds for the clique and independence numbers of graphs},
  author={Wilf, Herbert S},
  journal={Journal of Combinatorial Theory, Series B},
  volume={40},
  number={1},
  pages={113--117},
  year={1986},
  publisher={Elsevier},
url={https://doi.org/10.1016/0095-8956(86)90069-9}
}

@incollection{Smith1970,
  author    = {Smith, J. H.},
  title     = {Some properties of the spectrum of a graph},
  booktitle = {Combinatorial Structures and their Applications},
  publisher = {Gordon and Breach},
  address   = {New York},
  year      = {1970},
  pages     = {403--406}
}

@article{jadav2025strengthening,
  title={Strengthening Wilf's lower bound on clique number},
  author={Jadav, Hareshkumar and Madyastha, Sreekara and Raut, Rahul and Singh, Ranveer},
  journal={arXiv preprint},
  year={2025},
url={https://doi.org/10.48550/arXiv.2504.04836}
}

@article{Oboudi,
  title={Characterization of graphs with exactly two non-negative eigenvalues},
  author={Oboudi, Mohammad Reza},
  journal={Ars Mathematica Contemporanea},
  volume={12},
  number={2},
  pages={271--286},
  year={2016},
  url={https://doi.org/10.26493/1855-3974.1077.5b6}
}

@article{ando2015proof,
  title={Proof of a conjectured lower bound on the chromatic number of a graph},
  author={Ando, Tsuyoshi and Lin, Minghua},
  journal={Linear Algebra and its Applications},
  volume={485},
  pages={480--484},
  year={2015},
  publisher={Elsevier},
url={https://doi.org/10.1016/j.laa.2015.08.007}
}

@article{elphick2024two,
  title={Two conjectured strengthenings of Tur{\'a}n's theorem},
  author={Elphick, Clive and Linz, William and Wocjan, Pawel},
  journal={Linear Algebra and its Applications},
  volume={684},
  pages={23--36},
  year={2024},
  publisher={Elsevier},
url={https://doi.org/10.1016/j.laa.2023.12.010}
}

@article{stanley1987bound,
  title={A bound on the spectral radius of graphs with e edges},
  author={Stanley, Richard P},
  journal={Linear Algebra and its Applications},
  volume={87},
  pages={267--269},
  year={1987},
  publisher={Elsevier},
url={https://doi.org/10.1016/0024-3795(87)90172-8}
}

@article{wu2017upper,
  title={Upper bounds for the achromatic and coloring numbers of a graph},
  author={Wu, Baoyindureng and Elphick, Clive},
  journal={Discrete Applied Mathematics},
  volume={217},
  pages={375--380},
  year={2017},
  publisher={Elsevier},
url={https://doi.org/10.1016/j.dam.2016.09.005}
}

@article{favaron1993some,
  title={Some eigenvalue properties in graphs (conjectures of Graffiti—II)},
  author={Favaron, Odile and Mah{\'e}o, Maryvonne and Sacl{\'e}, J-F},
  journal={Discrete mathematics},
  volume={111},
  number={1-3},
  pages={197--220},
  year={1993},
  publisher={Elsevier},
url={https://doi.org/10.1016/0012-365X(93)90156-N}
}

@article{wocjan2012new,
  title={New spectral bounds on the chromatic number encompassing all eigenvalues of the adjacency matrix},
  author={Wocjan, Pawel and Elphick, Clive},
  journal={arXiv preprint},
  year={2012},
url={https://doi.org/10.48550/arXiv.1209.3190}
}

@article{edwards1983lower,
  title={Lower bounds for the clique and the chromatic numbers of a graph},
  author={Edwards, CS and Elphick, CH},
  journal={Discrete Applied Mathematics},
  volume={5},
  number={1},
  pages={51--64},
  year={1983},
  publisher={Elsevier},
url={https://doi.org/10.1016/0166-218X(83)90015-X}
}

@article{Duan23,
  title={A complete characterization of graphs with exactly two positive eigenvalues},
  author={Duan, Fang and Huang, Qiongxiang and Huang, Xueyi and Stani{\'c}, Zoran and Wang, Jianfeng},
  journal={Advances in Applied Mathematics},
  volume={144},
  pages={102457},
  year={2023},
  publisher={Elsevier},
  url={https://doi.org/10.1016/j.aam.2022.102457}
}

\end{document}